\documentclass[11pt]{article}
\usepackage{amsthm}
\usepackage{amsmath}
\usepackage{amssymb}
\usepackage{latexsym}
\usepackage{amsfonts}
\usepackage{mathrsfs}
\usepackage{listings}
\newtheorem{theorem}{\bf Theorem}[section]

\newtheorem{lemma}[theorem]{\bf Lemma}
\usepackage[none]{hyphenat}
\begin{document}
\title{$k$-Fibonacci numbers which are Narayana's cows numbers}
\author{Hayat Bensella and Djilali Behloul}
\date{}
\maketitle
\begin{abstract} \noindent
In this paper,we find all generalized Fibonacci numbers which are Narayana's cows numbers. In our proofs, we use both Baker's theory of nonzero linear forms in logarithms of algebraic numbers and the Baker-Davenport reduction method.
\end{abstract}

\section{Introduction}
\noindent
Let $\{F_n\}_{n\geq 0}$ be the Fibonacci sequence given by 
\[
F_{n+2} = F_{n+1} + F_n~~ {\rm for} ~n\geq 0
\]
with initials $F_0 = 0$ and $F_1 = 1$. 

Let $k \geq 2$ be an integer. Among the several generalizations of Fibonacci sequence, called the $k$-generalized Fibonacci sequence $\{F_{n}^{(k)} \}_{n \geq -(k-2)}$ we identify the one given by the recurrence

\begin{equation*}
  F_{n}^{(k)} = F_{n-1}^{(k)} + F_{n-2}^{(k)} + \dots + F_{n-k}^{(k)} = \sum_{i=1}^{k} F_{n-i}^{(k)} ~~\text{for all}~n \geq 2,  
\end{equation*}
with initial conditions $F_{-(k-2)}^{(k)} = F_{-(k-3)}^{(k)} = \dots = F_{0}^{(k)} = 0$ and $F_{1}^{(k)} = 1$. Here, $F_{n}^{(k)}$ denotes the $n$th $k$-generalized Fibonacci number.

Narayana's cows numbers $\{N_m\}_{n\ge 0}$ are intrduced by the indian mathematician Narayana Pandit In his book Ganita Kaumudi with the recurrence relation 
\[N_{m+3} = N_{m+2} + {N_m} \text{~~for~~} n\ge 0, \] 
 with initial terms $N_0 = 0$ and $ N_1 = N_2 = 1$. \\
 The Narayana sequence's defining relation is nearly identical to the well-known Fibonacci sequence, but with a delay that turns it into a third-order linear recurrence sequence. This is connected to Allouche and Johnson's "delayed morphisms"\cite{Allouche}.

Recently, diophantine equations related to  Narayana's cows sequence has attracted the interest of numerous authors. In \cite{repdigits} Bravo et al. find all repdigits in Narayana cow's sequence, in addition to this, they came to certain results on the existence of Mersenne prime numbers as well as numbers with distinct blocks of digits .
In \cite{Ray}, authors find all of the Narayana numbers that are sums of two base b repdigits. In \cite{fermat},  Bhoi and Ray found the Narayana cow's numbers, which are fermat numbers, and the same authors in  \cite{xcordinate}  searched for  $x$ coordinate of Pell equation which are Narayana sequence. Recently \cite{Hamtat} found all powers of two which are sums of two Narayana cow‘s numbers consedering the initial values $N_0=N_1=N_2=1$. The following theorem states his primary result
\begin{theorem} \label{hamtat}
    The Diophantine equation 
\begin{equation}\label{eq 3.13}
 N_m=2^l.    
\end{equation}
  has only the solutions:
  $(m,l)\in \{ (4,1),(6,2)\}$. Namely, $N_4=2$ and $N_6=4$.
\end{theorem}
Also several authors have worked on problems related to intersection of linear recurrent seqences and generalized Fibonacci sequence for examples see \cite{Bravo} \cite{Rihane 1} \cite{Rihane}.
In this paper, we study $k-$Fibonacci numbers which are Narayana's cows numbers, we found also a result about power of two in Narayana's cows sequences. More precisely, our main results is the following

\begin{theorem} \label{thm2}
The only nontrivial solutions of the Diophantine equation 

\begin{equation}\label{eq 1.3}
  F^{(k)}_n = N_m 
\end{equation}
on nonnegative integer $n, k,m$ with $ k\ge 2$, are\\
 \[F^{(2)}_{4}=N_5=3, ~~F^{(3)}_{6}=N_9=13 \text {~~and~~} F^{(k)}_{4}=N_6=4 \text {~~with~~} k\ge 3.\]
 
\end{theorem}
\section{Auxiliary results}
\noindent
\subsection{$k$-Generalized Fibonacci number}
First, we recall some facts and properties of the k-generalized Fibonacci sequence.
Note that the characteristic polynomial of the $k$-generalized Fibonacci sequence is
\[
\Psi_{k}(x) = x^k - x^{k-1} - \dots - x - 1.
\]
$\Psi_{k}(x)$ is irreducible over $\mathbb{Q}[x]$ and has just one root outside the unit circle. It is real and positive, so it satisfies $\alpha(k) > 1$. The other roots are strictly inside the unit circle. Throughout this paper, $\alpha := \alpha(k)$ denotes that single root, which is located between $2(1-2^{-k})$ and $2$ (see \cite{Miyazaki}).  In order to simplify notation, we will eliminates dependence on $k$ of $\alpha$.  

Dresden \cite{Dresden} gave a simplified Binet-like formula
for $F_{n}^{(k)}$:
\begin{equation*}
    F_{n}^{(k)} = \sum_{i=1}^{k} \frac{\alpha_{i}-1}{2+(k+1)(\alpha_{i}-2)} \alpha_{i}^{n-1},
\end{equation*}
where $\alpha = \alpha_{1}, \dots, \alpha_{k}$ are the roots of $\Psi_{k}(x)$. In addition, he demonstrated that the contribution of roots inside the unit circle to the right-hand side of \eqref{eq0} is very small. Specifically, he demonstrated that
\[
\Big| F_{n}^{(k)} - \frac{\alpha -1 }{2+(k+1)(\alpha - 2)} \alpha^{n-1} \Big| < \frac{1}{2}~~{\rm for~ all}~~n \geq 1.
\]
This allows us to write
\begin{equation} \label{eq0}
 F_{n}^{(k)} = f_k(\alpha) \alpha ^{n-1}+e_k(n). 
\end{equation}
where $\lvert  e_k(n) \vert <\frac{1}{2}$ for all $k\ge 2$ and $n\ge 1$.\\
The following inequality is proved by Bravo and Luca \cite[Lemma~1]{Bravo}.
\begin{lemma}\label{lm 2.1}
The inequality 
\[
\alpha^{n-2} \leq F_{n}^{(k)} \leq \alpha^{n-1}.
\]
holds for all $n \geq 1$.
\end{lemma}
One may immediately notice that the first $k + 1$ non-zero terms in $F_{n}^{(k)}$ are powers of $2$, namely
\[
F_{1}^{(k)}=1, F_{2}^{(k)}=1, F_{3}^{(k)}=2, F_{4}^{(k)}=4, \dots, F_{k+1}^{(k)} = 2^{k-1},
\]
while the next term in the above sequence is $F_{k+2}^{(k)} = 2^k -1$. Thus, we have that
\begin{equation} \label{eq 1.2}
 F_{n}^{(k)} = 2^{n-2}~\text{holds for all}~2 \leq n \leq k+1.   
\end{equation}
Also in \cite[pp. 542, 543]{Jhon} the authors proved that for all $n \ge k + 2$ we have
\begin{equation}\label{eq 2.4}
  F_n^{(k)}=2^{n-2}(1+ \zeta) \text{~~where~~} \zeta < \frac{k}{2}.  
\end{equation}

\subsection{The Narayana's cows sequences}
 The Binet's formula for Narayana's cows sequence is 
\begin{equation}\label{eq 2.5}
  N_m=a\lambda^m+b\beta ^m+c\gamma^m   
\end{equation}
for all $n\ge 0$\\
where $\lambda, \beta, \gamma $ are the roots of the characteristic equation $f(x)=x^3-x^2-1$ with   

\[
\gamma= \Bar{\beta},~~\lvert \beta \vert =\lvert \gamma \vert<1,~~ a=\frac{\lambda}{(\lambda-\beta)(\lambda-\gamma)},~~ b=\frac{\beta}{(\beta-\lambda)(\beta - \gamma)}, ~~c=\frac{\gamma}{(\gamma-\lambda)(\gamma-\beta)}.
\]
we can rewrite the formula \eqref{eq 2.5} as 
\begin{equation}\label{eq 2.6}
    N_m=C_\lambda \lambda^{m+2}+C_\beta \beta^{m+2}+C_\gamma \gamma ^{m+2} \text{~~~~for all~~~~}m\ge 0. 
\end{equation}

where
\[
C_x=\frac{1}{x^3+2}, ~ x\in \{ \lambda,\beta, \gamma\}.
\]
and $C_\lambda$ has the minimal polynomial $31x^3-31x^2+10x-1$ over $\mathbb{Z}$ and all the roots of this polynomial are strictly inside the unit circle. We can calculate 
\[
C^{-1}_\lambda \approx 5.1479 ~~and~~ C_\beta \approx 0.40751.
\]
Also, by induction one can prove that the $n^{th}$ Narayana number satisfies the following relation
\begin{equation}\label{eq 2.7}
    \lambda^{m-2} \le N_m \le \lambda^{m-1} \text{~~for all~~ } n\ge 1.
\end{equation}
The following result is proved by S$\acute{\text{a}}$nchez and Luca \cite[Lemma 7]{Sanchez}. 
\begin{lemma} \label{luca} 
If $ r \ge 1 , T> (4r^{2})^{r}$, and  $T> a/(\log a)^{r}$. Then  
 \[ 
 a <2^{r}T(\log T)^{r}.
 \]
\end{lemma}

In order to prove our main result, we use a Baker-type lower bound for a non-zero linear forms in logarithms of algebraic numbers a few times. Before presenting a result of Matveev \cite{Matveev} about the general lower bound for linear forms in logarithms, we recall some fundamental notations from algebraic number theory.

Let $\eta$ be an algebraic number of degree $d$ with minimal primitive polynomial 
\[
f(X):= a_0 X^d+a_1 X^{d-1}+ \cdots +a_d = a_0 \prod_{i=1}^{d}(X- \eta^{(i)}) \in \mathbb{Z}[X],
\]
where the $a_i$'s are relatively prime integers, $a_0 >0$, and the $\eta^{(i)}$'s are conjugates of $\eta$. Then
\begin{equation}\label{eq03}
h(\eta)=\frac{1}{d}\left(\log a_0+\sum_{i=1}^{d}\log\left(\max\{|\eta^{(i)}|,1\}\right)\right)
\end{equation}
is called the \emph{logarithmic height} of $\eta$.

With the established notations, Matveev (see  \cite{Matveev} or  \cite[Theorem~9.4]{Bugeaud1}), proved the ensuing result.

\begin{theorem}\label{th mtv}
Assume that $\gamma_1, \ldots, \gamma_t$ are positive real algebraic numbers in a real algebraic number field $\mathbb{K}$ of degree $D$, $b_1, \ldots, b_t$ are rational integers, and 
\[
\Lambda :=\eta_1^{b_1}\cdots\eta_t^{b_t}-1,
\]
is not zero. Then
\[
|\Lambda| \geq \exp\left(-1.4\cdot 30^{t+3}\cdot t^{4.5}\cdot D^2(1+\log D)(1+\log B)A_1\cdots A_t\right),
\]
where
\[
B\geq \max\{|b_1|,\ldots,|b_t|\},
\]
and
\[
A_i\geq \max\{Dh(\eta_i),|\log \eta_i|, 0.16\}, ~ \text{for all} ~ i=1,\ldots,t.
\]
\end{theorem}


Another result which will play an important role in our proof is due to Dujella and Peth\"{o} \cite[Lemma~5 (a)]{Dujella}. 

\begin{lemma} \label{lm daven}
Let $M$ be a positive integer, let $p/q$ be a convergent of the continued fraction of the irrational $\tau$ such that $q > 6M$, and let $A,B,\mu$ be some real numbers with $A>0$ and $B>1$. Let $\epsilon:=||\mu q||-M||\tau q||$, where $||\cdot||$ denotes the distance from the nearest integer. If $\epsilon >0$, then there exists no solution to the inequality
\[
0< |u \tau-v+\mu| <AB^{-u},
\]
in positive integers $u$ and $v$ with
\[
u \leq M \quad\text{and}\quad u \geq \frac{\log(Aq/\epsilon)}{\log B}.
\]
\end{lemma}

The following result is a simple property of the exponential function for further reference.
\begin{lemma}\label{lem exp}
For any non-zero real number $x$, we have the following

(a) $0 < x < |e^{x} - 1| $.

(b) If $x < 0$ and $|e^{x} - 1| < 1/2$, then $|x| < 2|e^{x} - 1|$.
\end{lemma}

\section{Proofs}
As an initial step, we consider the following observations.
\[
F^{(k)}_0=N_0,~ F^{(k)}_1=F^{(k)}_2=N_1=N_2=N_3=1,~ F^{(k)}_3=N_4=2.
\]
are valid for all $k\ge 2$. Hence the triples \[(n,k,m)\in \{(0,k,0),(1,k,1,(1,k,2),(1,k,3),(2,k,1),(2,k,2),(2,k,3),(3,k,4) \}\]
will be regarded as  trivial solutions of Equation \eqref{eq 1.3} for all $k\ge 2$, we assume throughout that $n \ge 4$, $m \ge 5$ and $k \ge 2$.

\subsection{The case $2\le n \le k+1$}
By the fact that  \eqref{eq 1.2} holds for $2\le n \le k+1$ we can rewrite \eqref{eq 1.3} as  
\begin{equation}\label{eq 3.12}
2^{n-2}=N_m .   
\end{equation}
Therefore the solution to our problem is to find powers of two in Narayana's cows sequence. Assume throughout that equation \eqref{eq 3.12} holds with $n\ge 4$ and $m \ge 5$. Suppose further that $k\ge 2$ in view of Theorem \ref{hamtat} that \eqref{eq 3.12} has only the solution $(n,k,m)=(4,k,6)$ with $k\ge 3$.
\subsection{The case $n\ge k+2$}
Assuming now that  $n\ge k+2$.  
Combining \eqref{eq 1.3}, \eqref{eq 2.7} and Lemma \ref{lm 2.1} to get 
\[
\alpha^{n-2}\le \lambda^{m-1} \text{~~and~~}\lambda^{m-2}\le \alpha^{n-1},
\]
which implies 
\[
(n-2)\frac{\log \alpha }{\log \lambda}+1 \le m \le (n-1)\frac{\log \alpha}{\log \lambda}+2.
\]
using $ \alpha >7/4$ for all $k\ge 2$, we obtain 
\begin{equation}\label{eq 3.20}
    1.4n-1.95< m< 1.9n+0.16 <2n.
\end{equation}
\subsubsection{An inequality for n and m in terms of k}
From \eqref{eq 1.3},\eqref{eq0} and \eqref{eq 2.6} we have 
\[
f_k(\alpha)\alpha^{n-1}-C_\lambda \lambda^{m+2}= C_\beta \beta ^{m+2}+C_\gamma \gamma^{m+2}-e_k(n).
\]
Taking absolute value on both sides 
\[
\lvert f_k(\alpha)\alpha^{n-1}-C_\lambda \lambda^{m+2} \vert < 2.
\]
Dividing  by $C_\lambda \lambda^{m+2}$
\begin{equation}\label{eq 3.21}
    \lvert C_\lambda ^{-1}f_k(\alpha)\alpha^{n-1}\lambda^{-(m+2)}-1\vert<\frac{2}{\lambda^m}.
\end{equation}
Put 
\[
\Lambda_2=C_\lambda ^{-1}f_k(\alpha)\alpha^{n-1}\lambda^{-(m+2)}-1. 
\]
In this application of Matveev's Theorem we take
\[
t:=3,~~ \eta_1= f_k(\alpha)/C_\lambda,~~\eta_2=\alpha,~~ \eta_3= \lambda,~~ b_1=1,~~ b_2=n-1,~~ b_3=-(m+2).
\]
Since $\eta_1,\eta_2,\eta_3 \in \mathbb{K}=\mathbb{Q}(\alpha,\lambda)$, we can take $D_{\mathbb{K}}=[\mathbb{K}:\mathbb{Q}]\ge 2k$, now we need to prove that $\Gamma_2 \neq 0$,  $f_k(\alpha)=C_\lambda \lambda ^{m+2}\alpha^{1-n}$, then $f_k(\alpha)$ would be an algebraic integer which is a contradiction (see \cite{Bravo}), thus $\Gamma_2 \neq 0$. The logarithmic heights are 
\[
h(\eta_1)=h(f_k(\alpha)/C_\lambda)\le h(f_k(\alpha))+h(C_\lambda)<\log(k+1)+\frac{\log 31}{3}<3.3 \log k.
\]
\[
h(\eta_2)=h(\alpha)= \frac{\log \alpha}{k}<\frac{\log 2}{k}.
\]
\[
h(\eta_3)=h(\lambda)=\frac{\log \lambda}{3}.
\]
We can take 
\[
A_1 \ge \max \{ 6.6k \log k,~~ \log (f_k(\alpha)/C_\lambda,~~0.16)  \}=6.6k\log k.
\]
\[
A_2 \ge \max \{ 2k \log 2/k,~~ \log (\alpha),~~0.16  \}=2\log 2.
\]
\[
A_3 \ge \max \{ 2k \log \lambda/3,~~ \log (\lambda),~~0.16  \}=\frac{2}{3}k \log \lambda<k\log \lambda.
\]
and  $B\ge \max\{1,n-1,m+2\}$,\eqref{eq 3.20} show that $m<1.9n+0.16<2n$ holds for $n\ge 5$, hence we can take $B:=2n+2$. Theorem \ref{th mtv} gives 
\[
\log \lvert \Lambda_2 \vert > -1.4\cdot 30^6\cdot 3^{4.5}\cdot(2k)^2(1+\log(2k))(1+\log(2n+2))6.6k\log k\cdot2\log2\cdot k\log\lambda.
\]
Comparing the above innequality with \eqref{eq 3.21} we get 
\[
m\log\lambda-\log2<2.02\cdot10^{12}k^4\log k (1+\log(2k))(1+\log(2n+2)).
\]
we use the fact $ 1+\log 2k < 3.5\log k $ and $1+\log (2n+2) <2.2\log n$ which hold for $k\ge2$ and $n\ge 5$, to get 
\begin{equation}\label{eq 3.22}
  m <3.94\cdot 10^{13}k^4 \log ^2k\log n. 
\end{equation}
By inequality \eqref{eq 3.20} we obtain 
\begin{equation*}\
    n<2.82\cdot 10^{13}K^4 \log ^2 k \log n,
\end{equation*}
which can be rewritten as 
\[
\frac{n}{\log n}< 2.82\cdot 10 ^{13} k^4\cdot \log ^2 k
\]
Applying Lemma \ref{luca} we get 
\begin{align*}
    n<& 2(2.82\cdot 10 ^{13} k^4\cdot \log ^2 k) \log (2.82\cdot 10 ^{13} k^4\cdot \log ^2 k)\\ < & 5.64 \cdot 10^{13} k^4 \log ^2 k(30.98+4\log k+ \log (\log k). 
\end{align*}
Using the fact that $30.98+4\log k+ \log (\log k)<49 \log k$, we get 
\begin{equation}\label{eq 3.23}
    n<2.77 \cdot 10^{15} k^4 \log ^3 k.
\end{equation}
\subsubsection{The case $2 \le k \le 220$}

Put 
\[
\Gamma_2=(n-1) \log \alpha -(m+2)\log \lambda + \log ( f_k(\alpha) C_\lambda).
\]
Then we can rewirte \eqref{eq 3.21} as 
\[
\lvert e^{\Gamma_1}-1\vert<\frac{2}{\lambda^m}.
\]
 $\Gamma_2 \neq 0$, since $\Lambda _2 \neq 0$, hence we discard the following cases:\\
If $\Gamma_1>0$, then $e^{ \Gamma_2}-1>0$, using Lemma \ref{lem exp}, we obtain 
\[
0<\Gamma_2<\frac{2}{\lambda^m}.
\]
If $\Gamma_2<0$, we have $\frac{2}{\lambda^m}< \frac{1}{2}$ holds for all $m \ge 5$, thus we get $\lvert e^{\Gamma_2}-1 \vert<\frac{1}{2} $ from Lemma \ref{lem exp}  again we have
\[
\Gamma_2< 2 \lvert e^{\Gamma_2}-1 \vert< \frac{4}{\lambda ^m} .
\]
In boths cases we get
\begin{equation} \label{eq 3.24}
    0< \Bigm \lvert (n-1) \frac{\log \alpha}{\log \lambda}-m + \frac{\log (f_k(\alpha)/C_\lambda \lambda^2)}{\log \lambda} \vert <5.2 \lambda ^{-m}.
\end{equation}
Puting 
\[
\tau_k=\frac{\log \alpha}{\log \lambda },~~ \mu_k=\frac{\log (f_k(\alpha)/C_\lambda \lambda^2)}{\log \lambda},~~ A=5.2,~~B=\lambda .
\]
It is clear that $\tau$ is an irrational number. We put $M_k=2.77\cdot 10^{15} k^4 \log ^3k$ which is an upper bound on $n-1$ from \eqref{eq 3.23}. Applying Lemma \ref{lm daven} for $k\in[2,220]$.\\
$\max M_k \approx 7\cdot 10^{26}$.  Let $q_{(m,k)}$ be the denominator of the $t$th convergent of the continued fraction of $\tau_{k}$.

We use {\it Mathematica} to get 
\[
\min_{\substack {2 \leq k \leq 220}} q_{(73,k)} > 10^{29} > 6M \mbox{ \ and \ } \max_{\substack {2 \leq k \leq 220}} q_{(73,k)} <  1.22 \cdot10^{43}. \
\]
The maximal value of $M\lVert \tau_{k} \cdot q_{(73,k)} \rVert < 0.0058$, whereas the minimal value of $\lVert \mu_{k} \cdot q_{(73,k)} \rVert> 0.0027$. 
Also, for 
\[
\epsilon_{73,k}:= \lVert\mu_{k} \cdot q_{(73,k)}\rVert - 7 \cdot 10^{26} \lVert\tau_{k} \cdot q_{(73,k)} \rVert,
\]
we obtain that 
\[
 \epsilon_{73,k} > 0.0027 , 
\]
 By Lemma \ref{lm daven}, we have for \eqref{eq 3.24}
\[
m \le \left \lfloor \frac{\log(5.2 \cdot1.22 \cdot 10^{43}/0.027}{\log 1.46}\right \rfloor .
\]
  therefore $m \leq 277$. From \eqref{eq 3.20} we get $n \leq 200$. Using that $k \leq n-2 $, we get $k \leq 198$. 
A computer search with {\it Mathematica}  in the following range:
\[ 
35\leq m\leq 277,  2\leq k \leq 198 ,  \text{ and }~ 4 \leq m \leq 200.
\]
we get th solutions listed in Theorem \ref{thm2}
This concludes the investigation of this case.
\subsubsection{The case $k>220$}

By using \eqref{eq 1.3},\eqref{eq 2.5} and \eqref{eq 2.4} we obtain 
\[
2^{n-2}-C_\lambda \lambda^{m+2}= C_\beta \beta ^{m+2}+C_\gamma \gamma^{m+2}- \zeta 2^{n-2}.
\]
Taking absolute value on both sides 
\[
\lvert2^{n-2}-C_\lambda \lambda^{m+2} \vert < \frac{1}{2}+ \frac{2^{n-2}}{2^{k/2}}.
\]
Dividing the above inequality by $2^{n-2}$
\begin{equation}\label{eq 3.25}
    \lvert 1-2^{-(n-2)}C_\lambda \lambda^{m+2}-1\vert<\frac{2}{2^{k/2}}.
\end{equation}
Put 
\[
\Lambda_3=1-2^{-(n-2)}C_\lambda \lambda^{m+2}-1.
\]
Appliying Matveev's Theorem with the data 
\[
t:=3,~~ \eta_1= 2,~~\eta_2=C_\lambda,~~ \eta_3= \lambda,~~ b_1=-(n-2),~~ b_2=1,~~ b_3=m+2.
\]
note that $\eta_1,\eta_2,\eta_3 \in \mathbb{K}=\mathbb{Q}(\lambda)$, we can take $D_{\mathbb{K}}=[\mathbb{K}:\mathbb{Q}]=3$, Now we show that $\Gamma_3 \neq 0$, indeed if this is were zero , we would then get $C_\lambda\lambda^{m+2}= 2^{n-2}  $ conjugating in $\mathbb{Q}(\lambda)$, we get $ \lvert C_\beta\beta^{m+2} \vert = 2^{n-2}$ leads to $\lvert C_\beta \vert >1 $ which is a contradiction thus $\Gamma_3 \neq 0$. The logarithmic heights are given by 
\[
h(\eta_1)=h(2)=\log 2.
\]
\[
h(\eta_2)=h(C_\lambda)= \frac{\log 31}{3}.
\]
\[
h(\eta_3)=h(\lambda)=\frac{\log \lambda}{3}.
\]
we can take 
\[
A_1 \ge \max \{ 3 \log 2,~~ \log 2,~~0.16)  \}=3\log 2.
\]
\[
A_2 \ge \max \{3\frac{\log 31}{3},~~ \log (C_\lambda),~~0.16  \}=3.44.
\]
\[
A_3 \ge \max \{ 3 \log \lambda/3,~~ \log (\lambda),~~0.16  \}=\frac{2}{3}k \log \lambda=0.39.
\]
and $B\ge \max\{n-2,1,m+2\}$, by \eqref{eq 3.20} we have  $m+2<1.9n+0.16+2<2n+2$ holds for $n\ge 5$ so we can take $B:=2n+2$. Therefore by Theorem \ref{th mtv} it result 
\begin{align*}
\log \lvert \Lambda_3 \vert >& -1.4\cdot 30^6\cdot 3^{4.5}\cdot(3)^2(1+\log3)(1+\log(2n+2))3.44 \cdot 0.39 \cdot 3 \log 2\\ >& -1.59\cdot 10^{13} \log n. 
\end{align*}
where we used the fact $1+\log (2n+2)<2.2 \log n$ holds for $n\ge 4$.\\
comparing the above innequality with \eqref{eq 3.25} we obtain 
\[
k<4.6\cdot10^{13}\log n .
\]
From \eqref{eq 3.23} we get 
\begin{align*}
  k<&4.6\cdot10^{13} \log( 2.77 \cdot 10^{15 k^4 \log^3k})\\<& 4.6\cdot10^{13}(35.6+4\log k+3 \log \log k)\\<&5.48\cdot 10^{14} \log k. 
\end{align*}
where we used $35.6+4\log k+3 \log \log k <11.9 \log k$. Using lemma \ref{luca} to get 
\[
k<3.72\cdot 10^{16}.
\]
From \eqref{eq 3.23} we get  
\begin{equation}\label{eq 3.26}
  n <2.96\cdot 10^{86} \text{~~and~~ }  m <5.92\cdot 10^{86}.
\end{equation}
Put 
\[
\Gamma_3=(m+2) \log \lambda -(n-2)\log 2 + \log  C_\lambda.
\]
Then \eqref{eq 3.25} can be rewirreten as 
\[
\lvert e^{\Gamma_3}-1\vert<\frac{2}{2^{k/2}}.
\]
Note that $\Gamma_3 \neq 0$, since $\Lambda _3 \neq 0$, so we destighish the following cases:\\
If $\Gamma_3>0$ the $e^{ \Gamma_3}-1>0$, using Lemma \ref{lem exp}, we obtain 
\[
0<\Gamma_3<\frac{2}{2^{k/2}}.
\]
If $\Gamma_3<0$, we have $\frac{2}{2^{k/2}}< \frac{1}{2}$ holds for all $k \ge 220$, gives $\lvert e^{\Gamma_3}-1 \vert<\frac{1}{2} $ from Lemma \ref{lem exp}  again we have
\[
\Gamma_3< 2 \lvert e^{\Gamma_3}-1 \vert< \frac{4}{2^{k/2}} .
\]
Hence in boths cases we get 
\begin{equation} \label{eq 3.27}
    0< \Bigm \lvert m \frac{\log \lambda}{\log 2}-n + \frac{2\log 2 \lambda+ \log C_\lambda}{\log 2} \vert <6 \cdot 2 ^{-k/2}.
\end{equation}
Puting 
\[
\tau:=\frac{\log \lambda}{\log 2},~~ \mu:=\frac{2\log 2 \lambda+ \log C_\lambda}{\log 2},~~ A:=6,~~B:=2 .
\]
clearly  $\tau$ is irrational number. We put $M=5.92\cdot 10^{86} $. to Apply Lemma \ref{lm daven} we let $q_{}$ be the denominator of the $t$th convergent of the continued fraction of $\tau$.
 we use {\it Mathematica} to get 
\[
q_{178} \approx 4.29\cdot  10^{87} > 6M.
\]
We have $M\lVert \tau \cdot q_{178} \rVert < 0.14$, whereas $\lVert \mu\cdot q_{178} \rVert> 0.20$. And 
\[
\epsilon_{178}:= \lVert\mu \cdot q_{178}\rVert - 5.92 \cdot 10^{86} \lVert\tau\cdot q_{178} \rVert,
\]
we obtain that 
\[
 \epsilon_{178} > 0.0.06 , 
\]
Hence by Lemma \ref{lm daven}, there are no integer solutions for \eqref{eq 3.27} when 
\[
\left \lfloor \frac{\log(6 \cdot q_{178}/0.06}{\log 2}\right \rfloor \leq \frac{k}{2} .
\]
  and therefore we have $k \leq 596$. Consequently from \eqref{eq 3.23} we get $n \leq 9.13 \cdot 10 ^{28}$ and $m < 1.83 \cdot 10^{29}$. \\
  Repeating Lemma \ref{lm daven} with $M:=1.83\cdot 10^{29}$ we get $q_{71}=10555900978374790722282223722863$, $\epsilon> 0.40-0.018>0.38$ thus 
 \[
\left \lfloor \frac{\log(6 \cdot q_{71}/0.38}{\log 2}\right \rfloor \leq \frac{k}{2} . 
\] 
which gives $k< 217$. This is a contradiction to  our assumption. Thus Theorem \ref{thm2} is proved.

\vspace{10mm} \noindent \footnotesize
\begin{minipage}[b]{90cm}
\large{USTHB, Faculty of Mathematics, \\ 
LATN Laboratory, BP 32, El Alia, 16111, \\ 
Bab Ezzouar, Algiers, Algeria.\\
Email: hbensella@usthb.dz}
\end{minipage}

\vspace{05mm} \noindent \footnotesize
\begin{minipage}[b]{90cm}
\large{USTHB, Faculty of computer sciences,\\
BP 32, El Alia, 16111 Bab Ezzouar, Algiers, Algeria.\\
Email: dbehloul@yahoo.fr}
\end{minipage}

\begin{thebibliography}{99}
 \bibitem{Allouche} J.P. Allouche and T. Johnson,  Narayana's cows and delayed morphisms, In \emph{articles of 3rd Computer Music Conference} JIM96, France 1996. 
 \bibitem{fermat} K. Bhoi and P.K. Ray,  Fermat numbers in Narayana's cows sequence. \emph{Integers} (2022).
\bibitem{xcordinate} Kisan Bhoi Prasanta Kumar Ray, On the $x$-Coordinates of Pell quations which are Narayana numbers. \emph{Integers} (2022). 
\bibitem{Bravo}  J.J. Bravo and F. Luca, Powers of two in generalized Fibonacci sequences, \emph{Rev. Colombiana Mat.}. \textbf{46} (2012), 67--79.
\bibitem{Jhon} Bravo, Jhon J., Carlos A. Gomez, and Jose L. Herrera. "On the intersection of k-Fibonacci and Pell numbers." Bulletin of the Korean Mathematical Society 56.2 (2019): 535-547.
\bibitem{repdigits} J.J Bravo, P. Das and S. Guzmàan, Repdigits in Narayana's cows sequence and their conse-
quences, J. Integer Seq. 23 (2020), Article 20.8.7.


\bibitem{Bugeaud1} Y. Bugeaud, M. Mignotte and S. Siksek, Classical and modular approaches to exponential Diophantine equations. I. Fibonacci and Lucas perfect powers, \emph{Ann. of Math. (2)} \textbf{163} (2006), 969--1018. 

\bibitem{Bugeaud2} Y. Bugeaud, F. Luca, M. Mignotte and S. Siksek, Fibonacci numbers at most one away from a perfect power, \emph{Elem. Math.} \textbf{63} (2008), 65--75. 

\bibitem{Dresden} G.P. Dresden and Z. Du, A simplified Binet formula for $k$-generalized Fibonacci numbers, \emph{J. Integer Seq.} {\bf 17}(4) (2014).

\bibitem{Dujella} A. Dujella and A. Peth\H o, A generalization of a theorem of Baker and Davenport, \emph{Quart. J. Math. Oxford Ser.} \textbf{49} (1998), 291--306.
\bibitem{Hamtat} A. Hamtat, An exponential Diophantine equation involving Narayana cow’s numbers, Utilitas Mathematica, 120, 37--48,(2023)

\bibitem{Koshy} T. Koshy, \emph{Fibonacci and Lucas numbers with Applications}, John Wiley \& Sons, 2011.

\bibitem{Laurent} M. Laurent, M. Mignotte and Y. Nesterenko, Formes lin\'eaires en deux logarithmes et d\'eterminants d'interpolation, (French) (Linear forms in two logarithms and interpolation determinants),  \emph{J. Number Theory} {\bf 55} (1995), 285--321.






\bibitem{Matveev} E.M. Matveev, An explicit lower bound for a homogeneous rational linear form in the logarithms of algebraic numbers, II, \emph{Izv. Ross. Akad. Nauk Ser. Mat.} \textbf{64} (2000), 125--180. Translation in \emph{Izv. Math.} \textbf{64} (2000), 1217--1269.
\bibitem{Miyazaki} T. Miyazaki, A. Wolfram, Solving generalized Fibonacci recurrences, \emph{Fibonacci Quart.} \textbf{36} (1998), 129--145.


\bibitem{Rihane 1}  S.E. Rihane, On k-Fibonacci balancing and k-Fibonacci
Lucas-balancing numbers. Carpathian Math. Publ. {\bf 13.1} (2021), 259-271.
 \bibitem{Rihane} S.E. Rihane and A. Togbé, k-Fibonacci numbers which are Padovan or Perrin numbers. Indian J Pure Appl Math (2022).


\bibitem{Sanchez} S. G. S$\acute{\text{a}}$nchez and F. Luca, Linear combinations of factorials and $S$-units in a binary recurrence sequence. Ann. Math. Qu$\acute{\text{e}}$. {\bf 38} (2014), 169--188.


\bibitem{Ray} R.P. Kumar, K. Bhoi, and B.K. Patel. "Narayana numbers as sums of two base b repdigits." Acta et Commentationes Universitatis Tartuensis de Mathematica 26.2 (2022): 183-192.










\end{thebibliography}
\end{document}